\documentclass[10pt,twoside]{article}
\usepackage{graphicx,amsmath,amsfonts,latexsym}
\usepackage{epic,here}
\usepackage{Latex-document}

\newcommand{\Figure}[2]{#2}   
\newtheorem{proposition} {Proposition}

\newcommand  {\C}     {\mathbb{C}}
\newcommand  {\T}     {\mathbb{T}}
\renewcommand  {\phi}   {\varphi}
\newcommand  {\RE} {{\rm Re }}

\newcommand  {\cd}{{\cdot}}
\newcommand  {\D}{{\partial}}

\title{\bf Adaptive Finite Element Methods  \vskip -2mm
for Partial Differential Equations\thanks{The author acknowledges
the support by the German Research Association (DFG) through SFB
359 `Reactive Flow, Diffusion and Transport'.}\vskip 6mm}

\author{R. Rannacher\vspace*{-0.5cm}\thanks{Institute of Applied Mathematics,
University of Heidelberg, Im Neuenheimer Feld 293/294, D-69120
Heidelberg, Germany. E-mail: rannacher@iwr.uni-heidelberg.de }}

\date{\vspace{-8mm}}

\begin{document}
\maketitle

\thispagestyle{first} \setcounter{page}{717}

\begin{abstract}\vskip 3mm
The numerical simulation of complex physical processes requires the use of
economical discrete models. This lecture presents a general paradigm of deriving
a posteriori error estimates for the Galerkin finite element approximation
of nonlinear problems. Employing duality techniques as used in optimal control
theory the error in the target quantities is estimated in terms of weighted `primal'
and `dual' residuals. On the basis of the resulting local
error indicators economical meshes can be constructed which are tailored to the
particular goal of the computation. The performance of this {\it Dual Weighted Residual
Method} is illustrated for a model situation in computational fluid mechanics:
the computation of the drag of a body in a viscous flow, the drag minimization by
boundary control and the investigation of the optimal solution's stability.
\vskip 4.5mm

\noindent {\bf 2000 Mathematics Subject Classification:}
65N30, 65N50, 65K10.

\noindent {\bf Keywords and Phrases:} Finite element method,
Adaptivity, Partial differential equations, Optimal control,
Eigenvalue problems.
\end{abstract}

\vskip 12mm

\section{Introduction} \label{section_1}\setzero
\vskip-5mm \hspace{5mm }

Suppose the goal of a simulation is the computation or optimization of a certain
quantity $\,J(u)\,$ from the solution $\,u\,$ of a continuous model with accuracy
$\,TOL\,$, by using the solution $\,u_h\,$ of a discrete model of dimension $\,N$,
\[
    {\cal A}(u) = 0, \qquad {\cal A}_h(u_h) = 0.
\]
Then, the goal of adaptivity is the optimal use of computing resources, i.e.,
minimum work for prescribed accuracy, or maximum accuracy for prescribed work.
In order to reach this goal, one uses {\it a posteriori} error estimates
\[
    |J(u)\!-\!J(u_h)| \,\approx\, \eta(u_h) := \sum_{K\in\T_h} \rho_K(u_h)\omega_K,
\]
in terms of the local residuals $\,\rho_K(u_h)\,$ of the computed solution
and weights $\,\omega_K\,$ obtained from the solution of a linearized
{\it dual problem}. In the following, we will describe a general optimal control
approach to such error estimates in Galerkin finite element methods.
For earlier work on adaptivity, we refer to the survey articles \cite{Verfurth96},
\cite{AinsworthOden97} and \cite{ErikssonEstepHansboJohnson95}. The contents
of this paper is based on material from \cite{BeckerRannacher95}, \cite{BeckerRannacher01}
and \cite{BangerthRannacher02}, where also references to other recent work can be
found.

\section{Paradigm of a posteriori error analysis} \label{section_2}
\setzero\vskip-5mm \hspace{5mm }

We develop a general approach to a posteriori error estimation for Galerkin
approximations of variational problems. The setting uses as little assumptions as
possible.
Let $\,X\,$ be some function space and $\,L(\cdot)\,$ a differentiable functional
on $\,X\,$. We are looking for stationary points of $\,L(\cdot)\,$
determined by
\[
    L'(x)(y) = 0 \qquad \forall y\in X,
\]
and their Galerkin approximation in finite dimensional subspaces $\,X_h\subset X\,$,
\[
    L'(x_h)(y_h) = 0 \qquad \forall\, y_h\in X_h.
\]
For this situation, we have the following general result:
\vspace{1mm}

\begin{proposition}\label{proposition_1}
There holds the {\it a posteriori} error representation
\begin{gather}\label{basic_error_id}
    L(x) - L(x_h) = \tfrac{1}{2} L'(x_h)(x\!-\!y_h) + R_h,
\end{gather}
for arbitrary $\,y_h\in X_h$. The remainder $\,R_h\,$ is cubic in $\,e := x\!-\!x_h$,
\[
        R_h := \tfrac{1}{2} \int_0^1 L'''(x_h\!+\!se)(e,e,e)
        \,s(s\!-\!1)\,\rm{ds}.
\]
\end{proposition}

{\noindent}{\bf Proof } We sketch the rather elementary proof.
First, we note that
\begin{align*}
    L(x)-L(x_h) &= \int_0^1 L'(x_h\!+\!se)(e)\,{\rm ds} \\
    &\hspace{7mm} - \tfrac{1}{2}\big\{ L'(x_h)(e) + L'(x)(e)\big\}
    + \tfrac{1}{2}L'(x_h)(e).
\end{align*}
Since $\,x_h\,$ is a stationary point,
\[
    L'(x_h)(e) = L'(x_h)(x\!-\!y_h) + L'(x_h)(y_h\!-\!x_h)
    = L'(x_h)(x\!-\!y_h), \quad y_h\in X_h.
\]
Finally, using the error representation of the trapezoidal rule,
\[
    \int_0^1 f(s)\,{\rm ds} - \tfrac{1}{2} \big\{
    f(0)+f(1)\big\} = {\textstyle \frac{1}{2}}\int_0^1 f''(s)s(s\!-\!1)\,{\rm ds},
\]
completes the proof. Notice that the derivation of the error
representation (\ref{basic_error_id}) does not assume the uniqueness of the
stationary points. But the a priori assumption $\;x_h\rightarrow x\;(h\!
\rightarrow\! 0)\,$ makes this result meaningful.


\section{Variational equations} \label{section_3}
\setzero\vskip-5mm \hspace{5mm }

We apply the result of Proposition \ref{proposition_1} to the Galerkin approximation
of {\it variational equations} posed in some function space $\,V$,
\begin{gather}
    a(u)(\psi) = 0 \quad\forall\psi\in V.
\end{gather}
Suppose that some functional output $\,J(u)\,$ of the solution $\,u\,$ is to be computed
using a Galerkin approximation in finite dimensional subspaces $\;V_h\subset V$,
\begin{gather}
    a(u_h)(\psi_h) = 0 \quad\forall\psi_h\in V _h.
\end{gather}
The goal is now to estimate the error $\,J(u)\!-\!J(u_h)\,$. To this end, we
employ a formal Euler-Lagrange approach to embed the present situation into the
general framework laid out above. Introducing a `dual' variable $\,z\,$
(`Lagrangian multiplier'), we define the Lagrangian functional
$\,{\cal L}(u,z) := J(u) - a(u)(z)\,$.
Then, stationary points $\;\{u,z\}\in V\!\times\! V\,$ of $\,{\cal L}(\cdot,\cdot) \,$
are determined by the system
\[
    {\cal L}'(u,z)(\phi,\psi) = \left\{  \begin{array}{l}
                   J'(u)(\phi) - a'(u)(\phi,z) \\
                   - a(u)(\psi)
                   \end{array}
            \right\} = 0 \quad \forall\{\phi,\psi\}.
\]
The corresponding Galerkin approximation determines $\{u_h,z_h\}\in V_h\!\times\! V_h$
by
\[
    {\cal L}'(u_h,z_h)(\phi_h,\psi_h) = \left\{  \begin{array}{l}
                   J'(u_h)(\phi_h) - a'(u_h)(\phi_h,z_h)
                   \rangle \\
                   - \langle a(u_h)(\psi_h)
                   \end{array}
            \right\} = 0 \quad \forall\{\phi_h,\psi_h\}.
\]

{\noindent}Set $\,x := \{u,z\},\,x_h := \{u_h,z_h\}$, and $\,L(x) := {\cal L}(u,z)\,$.
Then,
\[
    J(u)-J(u_h) \;=\; L(x) + a(u)(z) - L(x_h) - a(u_h)(z_h).
\]

\begin{proposition}\label{proposition_2}
With the `primal' and `dual' residuals
\begin{align*}
    \rho(u_h)(\cdot) &\,:=\, - a(u_h)(\cdot), \\
    \rho^{\ast}(z_h)(\cdot) &\,:=\, J'(u_h)(\cdot) \!-\! a'(u_h)(\cdot,z_h),
\end{align*}
there holds the error identity
\begin{gather}\label{error_id_1}
    J(u)-J(u_h) \;=\; \tfrac{1}{2}\rho(u_h)(z\!-\!\psi_h)\;+\;
    \tfrac{1}{2}\rho^\ast(z_h)(u\!-\!\phi_h) \;+\; {\cal R}_h,
\end{gather}
for arbitrary $\,\phi_h,\psi_h\in V_h\,$. The remainder $\,{\cal R}_h\,$
is cubic in the primal and dual errors $\,e^u:=u\!-\!u_h\,$ and  $\,e^z:=z\!-\!z_h\,$.
\end{proposition}

The evaluation of the error identity (\ref{error_id_1}) requires guesses for
primal and dual solutions $\,u\,$ and $\,z\,$ which are usually generated by
post-processing from the approximations $\,u_h\,$ and $\,z_h\,$,
respectively. The cubic remainder term $\,{\cal R}_h\,$ is neglected.
We emphasize that the solution of the dual problem takes only a
`linear work unit' compared to the solution of the generally nonlinear primal
problem.

\section{Optimal control problems} \label{section_4}
\setzero\vskip-5mm \hspace{5mm }

Next, we apply Proposition \ref{proposition_1} to the approximation of optimal
control problems. Let $\,V\,$ be the `state space' and $\,Q\,$ the `control space'
for the optimization problem
\begin{gather}
    J(u,q) \rightarrow \min!\qquad
    a(u)(\psi) + b(q,\psi) = 0 \quad\forall \psi\in V.
\end{gather}
Its Galerkin approximation uses subspaces $\,V_h\!\times\! Q_h\subset
V\!\times\! Q$ as follows:
\begin{gather}
    J(u_h,q_h) \rightarrow \min!\qquad
    a(u_h)(\psi_h) + b(q_h,\psi) = 0 \quad\forall\psi_h\in V_h.
\end{gather}
For embedding this situation into our general framework, we again employ the
Euler-Lagrange approach introducing the Lagrangian functional
$\,{\cal L}(u,q,z) := J(u,q) - A(u)(z) - B(q,z)\,$.
Corresponding stationary points $\;x := \{u,q,z\}\in X := V\!\times\!
Q\!\times\!V\,$ are determined by the system (`first-order optimality condition')
\begin{gather}
    \left\{ \begin{array}{l}
    J_u'(u,q)(\phi) - a'(u)(\phi,z) \\
    J_q'(u,q)(\chi) - b(\chi,z) \\
    - a(u)(\psi) - b(q,\psi)
   \end{array}
    \right\} = 0 \qquad \forall\{\phi,\chi,\psi\}.
\end{gather}
The Galerkin approximation detrmines $\,x_h := \{u_h,q_h,z_h\}\in
X_h := V_h\!\times\!Q_h\!\times V_h$ in finite dimensional subspace
$\,V_h\subset V,\,Q_h\subset Q\,$ by
\begin{gather}
    \left\{ \begin{array}{l}
    J_u'(u_h,q_h)(\phi_h) - a'(u_h)(\phi_h,z_h) \\
    J_q'(u_h,q_h)(\chi_h) - b(\chi_h,z_h) \\
    - a(u_h)(\psi_h) - b(q_h,\psi_h)
    \end{array}
    \right\} = 0 \qquad \forall\{\phi_h,\chi_h,\psi_h\}.
\end{gather}
For estimating the accuracy in this discretization, we propose to use the
natural `cost functional' of the optimization problem, i.e., to estimate
the error in terms of the difference $\,J(u,q)\!-\!J(u_h,q_h)\,$. Then, from
Proposition \ref{proposition_1}, we immediately obtain the following result:

\begin{proposition}\label{proposition_3}
With the `primal', `dual' and `control' residuals
\begin{align*}
    \rho^\ast(z_h)(\cdot) &\,:=\, J_u'(u_h,q_h)(\cdot) - a'(u_h)(\cdot,z_h), \\
    \rho^{q}(q_h)(\cdot) &\,:=\, J_q'(u_h,q_h)(\cdot) - b(\cdot,z_h), \\
    \rho(u_h)(\cdot) &\,:=\, - a(u_h)(\cdot) - b(q_h,\cdot),
\end{align*}
there holds the a posteriori error representation
\begin{align}
\begin{split}
    J(u,q)\!-\!J(u_h,q_h) \,&=\,
    \; \tfrac{1}{2} \rho^\ast(z_h)(u\!-\!\phi_h)
    + \tfrac{1}{2}\rho^{q}(q_h)(q\!-\!\chi_h) \\
    &\hspace{5mm}+ \tfrac{1}{2}\rho(u_h)(z\!-\!\psi_h) \;+\; {\cal R}_h,
\end{split}
\end{align}
for arbitrary $\,\phi_h,\,\psi_h\in V_h\,$ and $\,\chi_h\in Q_h$.
The remainder $\,{\cal R}_h\,$ is cubic in the errors $\,e^u:=u\!-\!u_h\,$,
$\,e^q := q\!-\!q_h\,$, $\,e^z := z\!-\!z_h\,$.
\end{proposition}

We note that error estimation in optimal control problems requires
only the use of available information from the computed solution
$\,\{u_h,q_h,z_h\}\,$, i.e., no extra dual problem has to be solve.
This is typical for a situation where the discretization error is measured with
respect to the `generating' functional of the problem, i.e. the Lagrange functional
in this case.
In the practical solution process the mesh adaptation is nested with an outer
Newton iteration leading to a successive `model enrichment'.
The `optimal' solution $\,\{u_h^{\rm opt},q_h^{\rm opt}\}\,$ obtained by the adapted
discretization may satisfy the state equation only in a rather week sense.
If more `admissibility' is required, we may solve just the state equation with
an better discretization (say on a finer mesh) using the computed optimal control
$\,q_h^{opt}\,$ as data.

\section{Eigenvalue problems} \label{section_5}
\setzero\vskip-5mm \hspace{5mm }

Finally, we apply Proposition \ref{proposition_1} to the Galerkin approximation
of eigenvalue problems. Consider in a (complex) function space $\,V\,$
the generalized eigenvalue problem
\begin{gather}\label{eigenvalueproblem}
    a(u,\psi) = \lambda\, m(u,\psi) \quad\forall\psi\in V, \qquad \lambda\in \C, \;
    m(u,u) = 1,
\end{gather}
where the form $\,a(\cdot,\cdot)\,$ is linear but not necessarily symmetric, and the
eigenvalue form $\,m(\cdot,\cdot)\,$ is symmetric and positive semi-definit.
The Galerkin approximation is defined in finite dimensional subspaces $V_h\subset V$,
\begin{gather}\label{discr_eigenvalueproblem}
    a(u_h,\psi_h) = \lambda_h m(u_h,\psi_h)\quad\forall\psi_h\in V_h, \qquad
    \lambda_h\in \C,\, m(u_h,u_h) = 1.
\end{gather}
We want to control the error in the eigenvalues $\,\lambda\!-\!\lambda_h\,$. To this
end, we embed this situation into the general framework of variational equations
by introducing the spaces $\,{\cal V} := V\!\times\!\C\,$ and
$\,{\cal V}_h := V_h\!\times\!\C\,$, consisting of elements $\,U := \{u,\lambda\}\,$
and $\,U_h := \{u_h,\lambda_h\}\,$, and the semi-linear form
\[
   A(U)(\Psi) := \lambda m(u,\psi) \!-\! a(u,\psi)
    + \overline{\mu}\big\{ m(u,u) - 1\big\},
    \quad  \Psi=\{\psi,\mu\}\in {\cal V}.
\]
Then, the eigenvalue problem (\ref{eigenvalueproblem}) and its Galerkin approximation
(\ref{discr_eigenvalueproblem}) can be written in the compact form
\begin{eqnarray}
    A(U)(\Psi) &=& 0 \quad\quad\forall \Psi\in {\cal V}, \\
    A(U_h)(\Psi_h) &=& 0 \qquad\forall \Psi_h\in {\cal V}_h.
\end{eqnarray}
The error in this approximation will be estimated with respect to the functional
\[
    J(\Phi) := \mu\,m(\phi,\phi),
\]
where $\,J(U) = \lambda\,$ since $\,m(u,u) = 1\,$.
The corresponding continuous and discrete dual solutions $\;Z = \{z,\pi\}\in
{\cal V}\,$ and $\,Z_h = \{z_h,\pi_h\}\in {\cal V}_h\,$ are determined by the
problems
\begin{align}
    A'(U)(\Phi,Z) &= J'(U)(\Phi) \qquad\;\;\;\forall \Phi\in {\cal V}, \\
    A'(U_h)(\Phi_h,Z_h) &= J'(U_h)(\Phi_h) \qquad\forall \Phi_h\in {\cal V}_h.
\end{align}
A straightforward calculation shows that these dual problems are equivalent
to the {\it adjoint eigenvalue problems} associated to (\ref{eigenvalueproblem}) and
(\ref{discr_eigenvalueproblem}),
\begin{align}
    a(\phi,z) &= \pi\,m(\phi,z) \hspace{11mm} \forall \phi\in V,
    \hspace{8mm} m(u,z) = 1, \\
    a(\phi_h,z_h) &= \pi_h\,m(\phi_h,z_h) \quad \forall \phi_h\in V_h,
    \quad m(u_h,z_h) = 1.
\end{align}
Then, application of Proposition \ref{proposition_1} yields the following result:
%
\begin{proposition}\label{propopsition_4}
With the `primal' and `dual' residuals
\begin{align*}
    \rho(u_h,\lambda_h)(\cdot) &:= a(u_h,\cdot) \!-\! \lambda_h\,m(u_h,\cdot), \\
    \rho^\ast(z_h,\pi_h)(\cdot) &:= a(\cdot,z_h) \!-\! \pi_h\,m(\cdot,z_h),
\end{align*}
there holds the a posteriori error representation
\begin{gather}\label{eigenvalue_error}
    \lambda\!-\!\lambda_h
    \,=\,\tfrac{1}{2} \rho(u_h,\lambda_h)(z\!-\!\psi_h)
     + \tfrac{1}{2}\rho^\ast(z_h,\pi_h)(u\!-\!\phi_h) - {\cal R}_h ,
\end{gather}
for arbitrary $\,\psi_h,\,\phi_h\in V_h$, with the remainder term
\[
    {\cal R}_h \,=\, \textstyle{\frac{1}{2}}(\lambda\!-\!\lambda_h)\,
    m(v\!-\!v_h,z\!-\!z_h) .
\]
\end{proposition}

We note that in Proposition \ref{propopsition_4}, no assumption about the multiplicity
of the approximated eigenvalue $\,\lambda\,$ has been made. In order to make the
error representation (\ref{eigenvalue_error}) meaningful, we have to use a priori
information about the convergence $\,\{\lambda_h,v_h\}\rightarrow \{\lambda,v\}\,$ as
$\,h\rightarrow 0\,$. The simultaneous solution of primal and dual eigenvalue problems
naturally occurs within an optimal multigrid solver of nonsymmetric eigenvalue problems.
Further, error estimates with respect to functionals $\,J(u)\,$ of eigenfunctions can be
derived following the general paradigm. Finally, in solving {\it stability eigenvalue
problems} $\,{\cal A}'(\hat{u})v = \lambda{\cal M}v\,$, we can include the perturbation
of the operator $\,{\cal A}'(\hat{u}_h) \approx {\cal A}'(\hat{u})\,$ in the
a posteriori error estimate of the eigenvalues.

\section{Application in fluid flow simulation} \label{section_6}
\setzero\vskip-5mm \hspace{5mm }

In order to illustrate the abstract theory developed so far, we present some
results for the application of `residual-driven' mesh adaptation for a model
problem in computational fluid mechanics, namely  `channel flow around a cylinder'
as shown in the figure below. The {\it stationary} Navier-Stokes system
\[
    {\cal A}(u) :=  \left\{ \begin{array}{ll}
                    - \nu\Delta v + v\cd\nabla v + \nabla p\\
                    \hspace{20mm} \nabla\cd v
                    \end{array} \right\} = 0
\]
determines the pair $\,u := \{v,p\}\,$ of velocity vector $\,v\,$ and scalar pressure
$\,p\,$ of a viscous incompressible fluid with viscosity $\,\nu\,$ and normalized
density $\,\rho \equiv 1\,$. The physical boundary conditions are
$\,v|_{\Gamma_{\rm{rigid}}}=0\,$, $\,v|_{\Gamma_{\rm{in}}} = v^{\rm{in}}\,$, and
$\,\nu\partial_n v-np|_{\Gamma_{\rm{out}}}=0\,$,
i.e., the flow is driven by the prescribed parabolic inflow $\,v^{\rm{in}}\,$.
The Reynolds number is $\,\rm{Re} = \frac{\bar{U}^2D}{\nu} = 20\,$,
such that the flow is stationary.
%
\begin{figure}[H]
\begin{center}
\setlength{\unitlength}{900sp}
\begingroup\makeatletter\ifx\SetFigFont\undefined%
\gdef\SetFigFont#1#2#3#4#5{%
  \reset@font\fontsize{#1}{#2pt}%
  \fontfamily{#3}\fontseries{#4}\fontshape{#5}%
  \selectfont}%
\fi\endgroup%
\begin{picture}(15000,3772)(1,-4197)
\thinlines
\put(3151,-2311){\circle{900}}
\put(901,-511){\line( 0,-1){900}}
\put(901,-1411){\line( 0,-1){2700}}
\put(901,-4111){\line( 1, 0){12000}}
\put(901,-511){\line( 1, 0){12000}}
\put(901,-511){\makebox(1.5875,11.1125){\SetFigFont{5}{6}
{\rmdefault}{\mddefault}{\updefault}.}}
\thicklines
\put(2951,-511){\line( 1, 0){900}}
\put(2951,-4111){\line( 1, 0){900}}
\put(2951,-505){\line( 1, 0){900}}
\put(2951,-4109){\line( 1, 0){900}}
\put(2951,-500){\line( 1, 0){900}}
\put(2951,-4122){\line( 1, 0){900}}
\put(  1,-2761){\makebox(0,0)[lb]{\smash{\SetFigFont{12}{14.4}
{\rmdefault}{\mddefault}{\updefault}\hspace{-2mm}\large$\Gamma_{\rm{in}}$}}}
\put(3601,-2761){\makebox(0,0)[lb]{\smash{\SetFigFont{12}{14.4}
{\rmdefault}{\mddefault}{\updefault}\hspace{1mm}\large S}}}

\put(13000,-2761){\makebox(0,0)[lb]{\smash{\SetFigFont{12}{14.4}
{\rmdefault}{\mddefault}{\updefault}\hspace{1mm}\large$\Gamma_{\rm{out}}$}}}

\put(2251,-1150){\makebox(0,0)[lb]{\smash{\SetFigFont{12}{14.4}
{\rmdefault}{\mddefault}{\updefault}\large$\Gamma_{1}$}}}
\put(2251,-3861){\makebox(0,0)[lb]{\smash{\SetFigFont{12}{14.4}
{\rmdefault}{\mddefault}{\updefault}\large$\Gamma_{2}$}}}
\end{picture}
\end{center}
\end{figure}
%
Let the goal of the simulation be the accurate computation
of the effective force in the main flow direction imposed on the cylinder, i.e.
the so-called `drag coefficient',
\[
    J(u) := c_{\rm drag} = \frac{2}{\max|v^{\rm in}|^2 D} \int_S n^T(2\nu\tau \!-\!
    pI)e_1\,{\rm ds},
\]
where $\,S\,$ is the surface of the cylinder, $\,D\,$ its diameter, and
$\,\tau = \frac{1}{2}(\nabla v\!+\!\nabla v^T)\,$ the strain tensor.
In practice, one uses a volume-oriented representation of $\,c_{\rm drag}\,$.

Here, we cannot describe the standard variational formulation of the Navier-Stokes
problem and its Galerkin finite element discretization in detail but rather refer
to the literature; see \cite{Rannacher00}, \cite{BeckerRannacher01},
and the references therein.

In the present situation the primal and dual residuals occuring in the a posteriori
error representation (\ref{error_id_1}) have the following explicit form:
\begin{align*}
    &\rho(u_h)(z\!-\!z_h) := \sum_{K\in\T_h} \Big\{ (R_h,
    z^v\!-\!z_h^v)_K + (r_h,z^v\!-\!z_h^v)_{\D K}
    + (z^p\!-\!z_h^p,\nabla\cd v_h)_K + \dots \Big\},\\
    &\rho^\ast(z_h)(u\!-\!u_h) := \sum_{K\in\T_h} \Big\{ (R^\ast_h,v\!-\!v_h)_K
    + (r^\ast_h,v\!-\!v_h)_{\D K}
    + (p\!-\!p_h,\nabla\cd z_h^v)_K + \dots \Big\},
\end{align*}
with the cell and edge residuals defined by
\begin{align*}
    R_{h|K} &:= f+\nu\Delta v_h \!-\!v_h\cd\nabla v_h\!-\!\nabla p, \\
    R^\ast_{h|K} &:= j+\nu\Delta z^v_h \!+\! v_h\cd\nabla z^v_h \!-\!
    \nabla v_h^Tz_h^v \!+\! \nabla\cd v_hz_h^v \!-\!\nabla z^p_h,\\
    r_{h|\Gamma} &:=
    \left\{ \begin{array}{ll}
    \tfrac{1}{2}[\nu\D_n v_h\!-\!np_h], &\rm{if}\;\;\Gamma\not\subset
    \D\Omega \\
    - \nu\D_n v_h\!+\!np_h, &\rm{if}\;\;\Gamma\subset \Gamma_{\rm{out}},
    \quad (= 0\;\;\text{else})
    \end{array} \right\},\\
    r^\ast_{h|\Gamma} &:=
    \left\{ \begin{array}{ll}
    \tfrac{1}{2}[\nu\D_n z^v_h\!+\!n\cd v_hz_h^v\!-\!z^p_hn],
    &\rm{if}\;\;\Gamma\not\subset\partial\Omega \\
    - \nu\D_n z^v_h\!-\!n\cd v_hz_h^v\!+\!z^p_hn,
    &\rm{if}\;\;\Gamma\subset \Gamma_{\rm{out}},\quad (= 0\;\;\text{else})
    \end{array} \right\},
\end{align*}
where $\,[\dots]\;$ denots the jump across edges $\,\Gamma\,$, and
`$\dots$' stands for terms representing errors due to boundary and inflow
approximation as well as stabilization.

Practical mesh adaptation on the basis of the a posteriori error estimates
proceeds as follows:
At first, the error functional may have to be regularized according
to $\,\tilde{J}(u) = J(u) + {\cal O}(TOL)\,$.
Then, after having computed the primal approximation $\,u_h\,$, the {\it linear}
discrete dual problem is solved:
\begin{gather}
    \langle {\cal A}'(u_h)^{\ast}z_h,\phi_h \rangle = \tilde{J}'(u_h)(\phi_h) \quad
    \forall \phi_h\in V_h^\ast.
\end{gather}
The error estimator is localized, $\,\eta_\omega = \sum_{K\in\T_h}\eta_K\,$, and
approximation of the weights are computed by patch-wise higher-order interpolation:
$\,(z\!-\!z_h)_{|K} \,\approx\, (I_{2h}^\ast z_h\!-\!z_h)_{|K}\,$.
Finally, the current mesh is adapted by `error balancing'
$\,\eta_K \,\approx\, \eta_\omega/\#\{K\in\T_h\}\,$.
In the following, we show some results which have been obtained using mesh
adaptation on the basis of the Dual Weighted Residual Method (`DWR method').


\subsection{Drag computation (from \cite{Becker00})}
\setzero\vskip-5mm \hspace{5mm }

The drag is computed on meshes generated by the DWR method and by
an `ad hoc' refinement criterion based on smoothness properties
of the computed solution. \vspace{-10mm}

\begin{table}[H]
\begin{center}
\caption{\label{tab_11:drag-lift_results_1} Results for drag
computation on adapted meshes ($1\%$-error in bold face).}
\vspace{3mm}

\begin{tabular}{|c|c|c|c|c|} \hline
\multicolumn{5}{|c|}{Computation of drag} \\\hline
$L $&$      N  $&$ c_{\rm{drag}} $&$ \eta_{\rm{drag}} $&$ I_{\rm{eff}}$
\\\hline
$4 $&$    984  $&$  5.66058       $&$ 1.1e\!-\!1   $&$0.76          $ \\\hline
$5 $&${\bf 2244}$&$  5.59431       $&$ 3.1e\!-\!2   $&$0.47          $ \\\hline
$6 $&$   4368  $&$ 5.58980 $&$ 1.8e\!-\!2   $&$0.58          $ \\\hline
$6 $&$   7680  $&$  5.58507       $&$ 8.0e\!-\!3   $&$0.69          $ \\\hline
    &$  \infty $&$  5.57953       $&                &                 \\\hline
\end{tabular}
\end{center}
\end{table}
\vspace{-5mm}

\begin{figure}[H]
\begin{center}
\includegraphics[width=0.6\textwidth]{\Figure{figures}{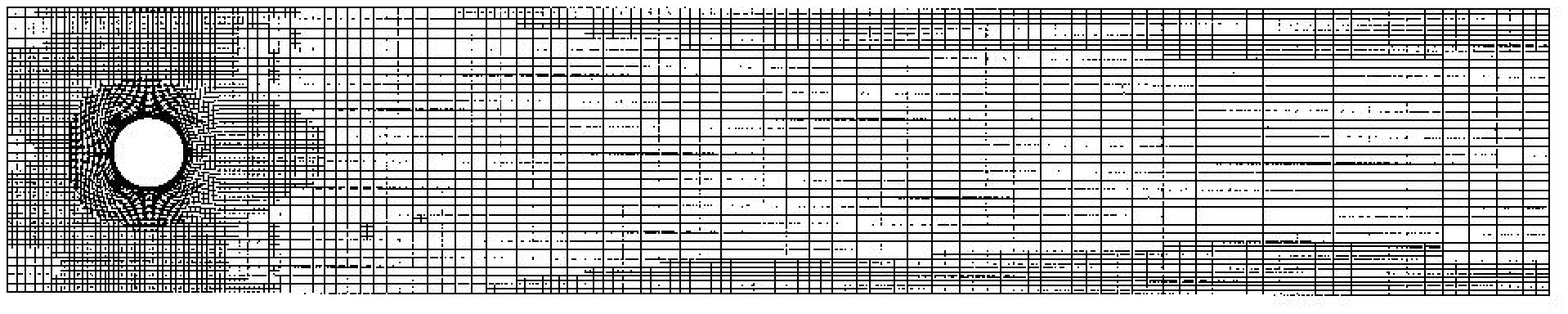}}
\includegraphics[width=0.6\textwidth]{\Figure{figures}{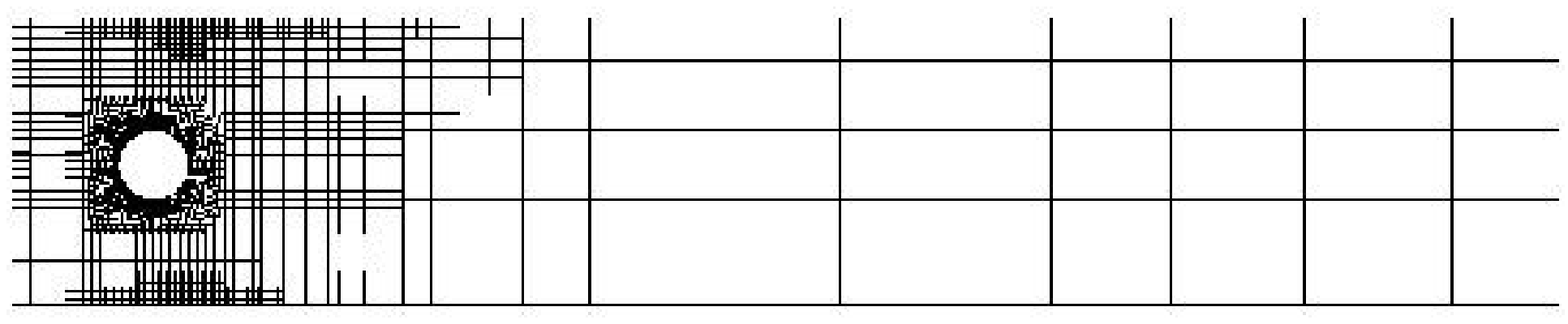}}
\begin{minipage}{10cm}
\caption{\label{fig_11:samplemeshes} Refined meshes by `ad hoc' strategy (top) and DWR method (bottom)}
\end{minipage}
\end{center}
\end{figure}
\vspace{-5mm}


\subsection{Drag minimization (from \cite{Becker01})}
\setzero\vskip-5mm \hspace{5mm }

The drag coefficient is to be minimized by imposing a pressure drop at the
two outlets $\,\Gamma_i\,$ above and below the cylinder. In this case of
`boundary control' the control form is given by
$\,b(q,\psi) := - (q,n\cd \psi^v)_{\Gamma_1\cup\Gamma_2}\,$.
\vspace{-7mm}

\begin{table}[H]
\begin{center}
\caption{\label{tab_11:drag-lift_results_2} Uniform refinement
versus adaptive refinement for $\,{\rm Re} = 40\,$.} \vspace{3mm}

\begin{tabular}{|c|c||c|c|} \hline
\multicolumn{2}{|c||}{Uniform refinement}&\multicolumn{2}{c|}
{Adaptive refinement} \\\hline
$N$ & $J_{\rm{drag}}$   & $N$ & $J_{\rm{drag}}$ \\\hline
$10512  $&$ 3.31321  $&$  1572$&$  3.28625 $ \\\hline
$41504  $&$ 3.21096  $&$  4264$&$  3.16723 $ \\\hline
$164928 $&$ 3.11800  $&$ 11146$&$  3.11972 $ \\\hline
\end{tabular}
\end{center}
\end{table}
\vspace{-10mm}

\begin{figure}[H]
\begin{center}
\rotatebox{-90}{
\includegraphics[scale=0.31]{\Figure{figures}{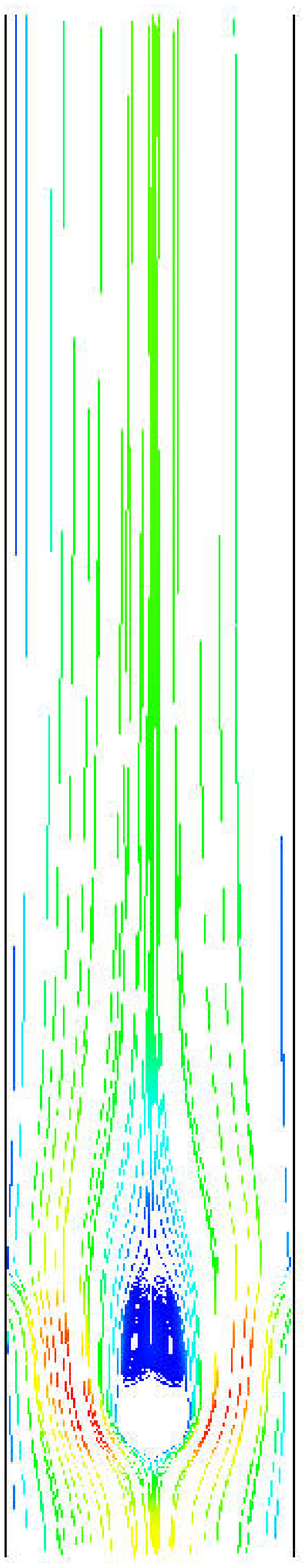}}
\hspace{-3mm}
\includegraphics[scale=0.31]{\Figure{figures}{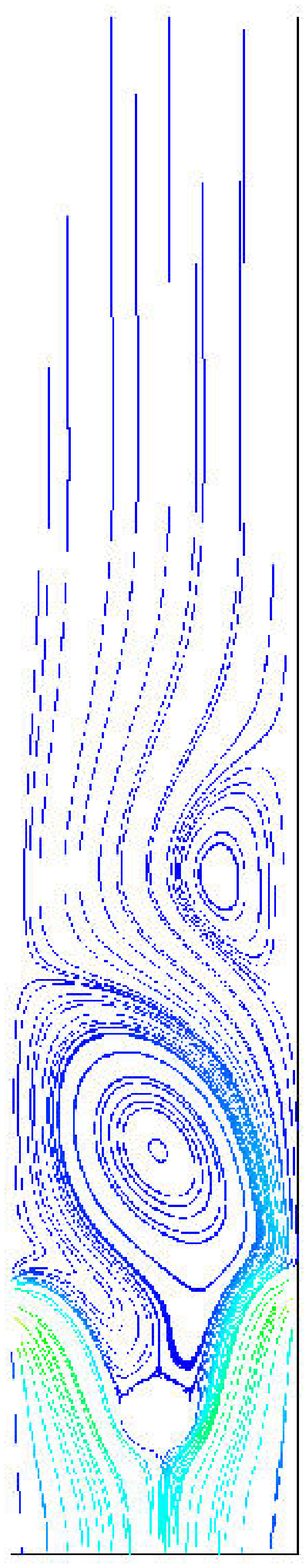}}
\hspace{-6mm}
\includegraphics[scale=0.32]{\Figure{figures}{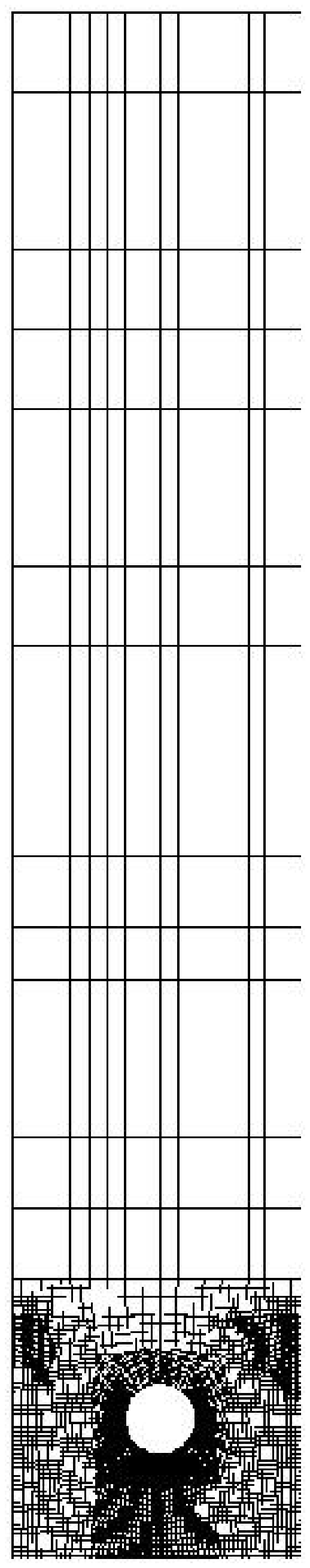}}
} \vspace{-5mm}
\begin{minipage}{10cm}
\caption{ Velocity of the uncontrolled flow (top), controlled flow (middle), corresponding adapted mesh (bottom)}
\end{minipage}
\end{center}
\end{figure}
\vspace{-5mm}


\subsection{Stability of optimized flows (from \cite{HeuvelineRannacher02})}
\vskip-5mm \hspace{5mm }

We want to investigate the stability of the optimized solution
$\,u^{\rm opt} = \{v^{\rm opt},p^{\rm opt}\}\,$ by linear stability theory.
This is a crucial question since in the present case the optimal
solution is obtained by a {\it stationary} Newton iteration which may converge
to physically unstable solutions. In this context, we have to consider the
non-symmetric eigenvalue problem for $\;u := \{v,p\}\in V\,$ and
$\,\lambda\in \C\,$:
\[
    {\cal A}'(u^{\rm opt})u :=  \left\{ \begin{array}{ll}
    -\nu\Delta v + v^{\rm opt}\cd\nabla v + v\cd\nabla v^{\rm opt} + \nabla p\\
    \hspace{20mm} \nabla\cd v
    \end{array} \right\}
    = \lambda \left\{ \begin{array}{l}
    v \\ 0
    \end{array} \right\}.
\]
If the real parts of all eigenvalues are positive, $\,\RE\,\lambda > 0\,$,
then the (stationary) base flow $\,\{v^{\rm opt},p^{\rm opt}\}\,$ is considered
as stable (but with respect to possibly only very small perturbations). We find
that the optimal solution is at the edge of being unstable.
\vspace{-2mm}

\begin{figure}[H]
\begin{center}
    \includegraphics[scale=0.25]{\Figure{figures}{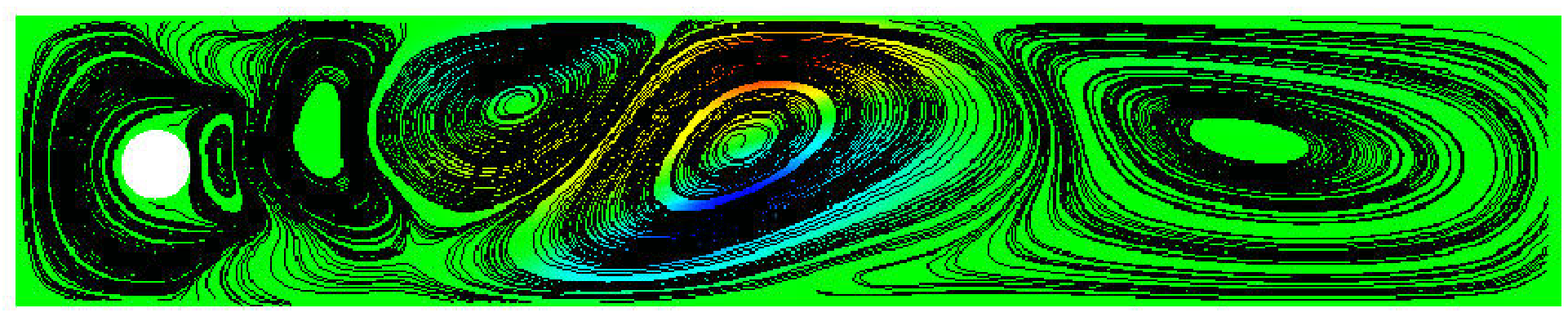}}
\vspace{3mm}

    \includegraphics[scale=0.25]{\Figure{figures}{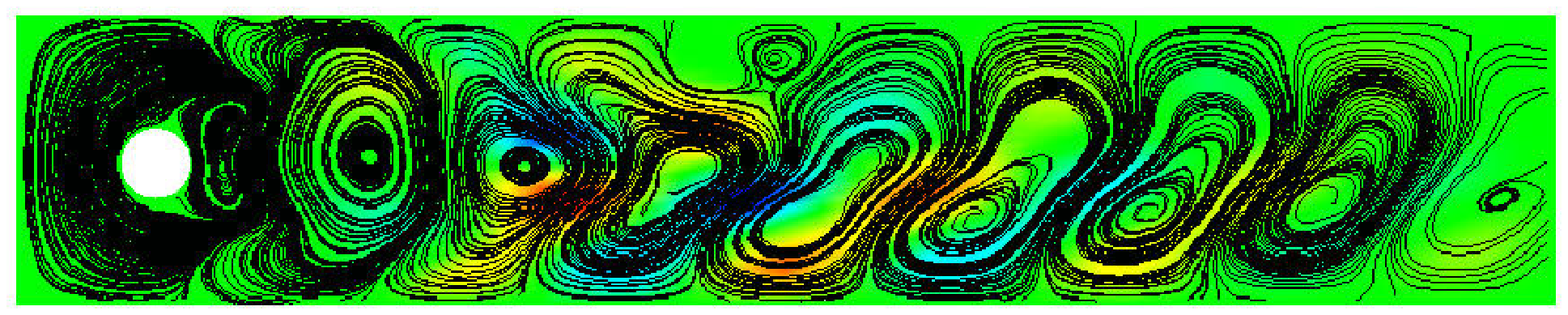}}
\vspace{-3mm}

\begin{minipage}{10cm}
\caption{\it Streamlines of real parts of the `critical' eigenfunction shortly before the Hopf bifurcation and
after, depending on the imposed pressure drop}
\end{minipage}
\end{center}
\end{figure}
\vspace{-5mm}

\begin{figure}[H]
\begin{center}
    \includegraphics[scale=0.30]{\Figure{figures}{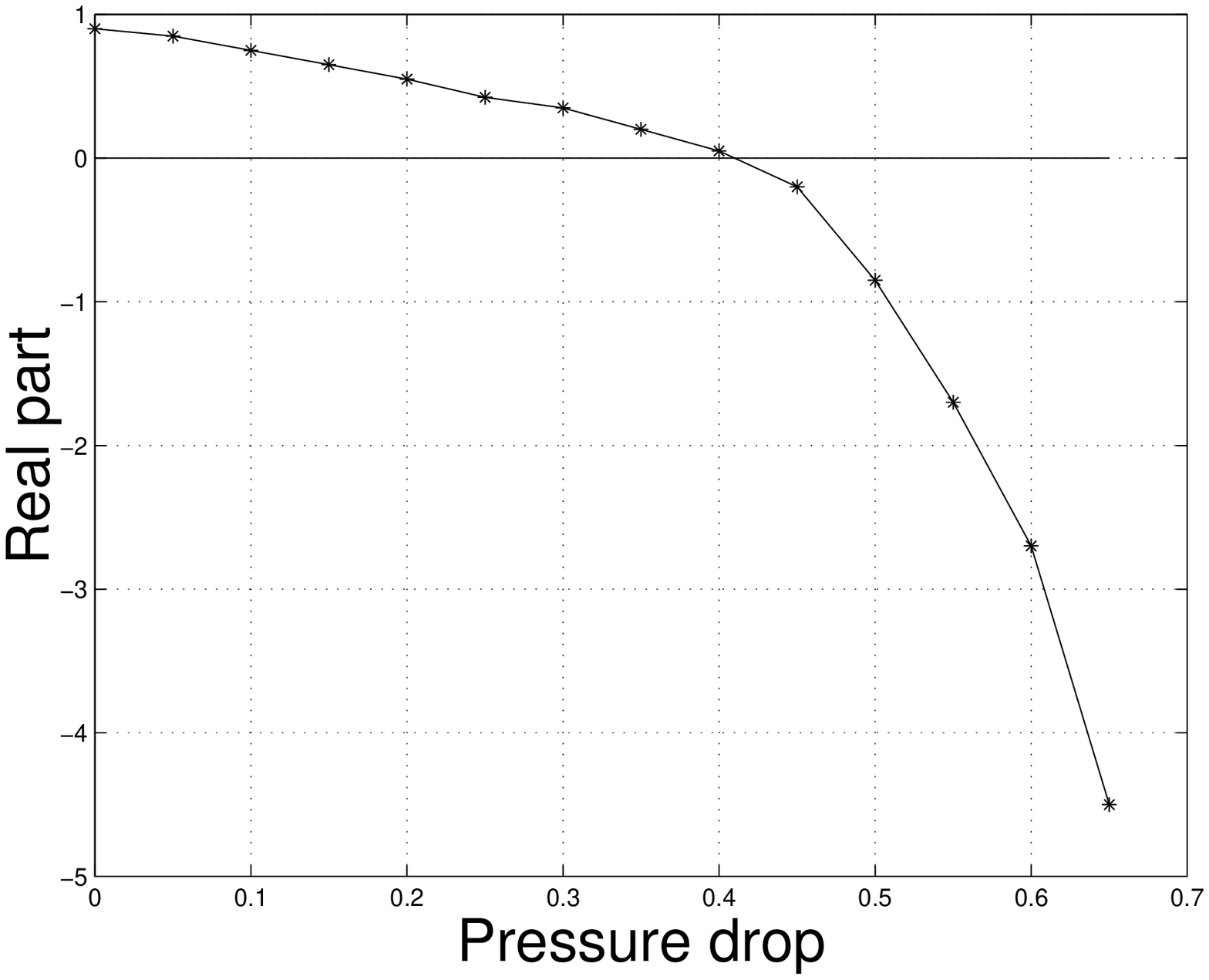}}
\hspace{10mm}
    \includegraphics[scale=0.30]{\Figure{figures}{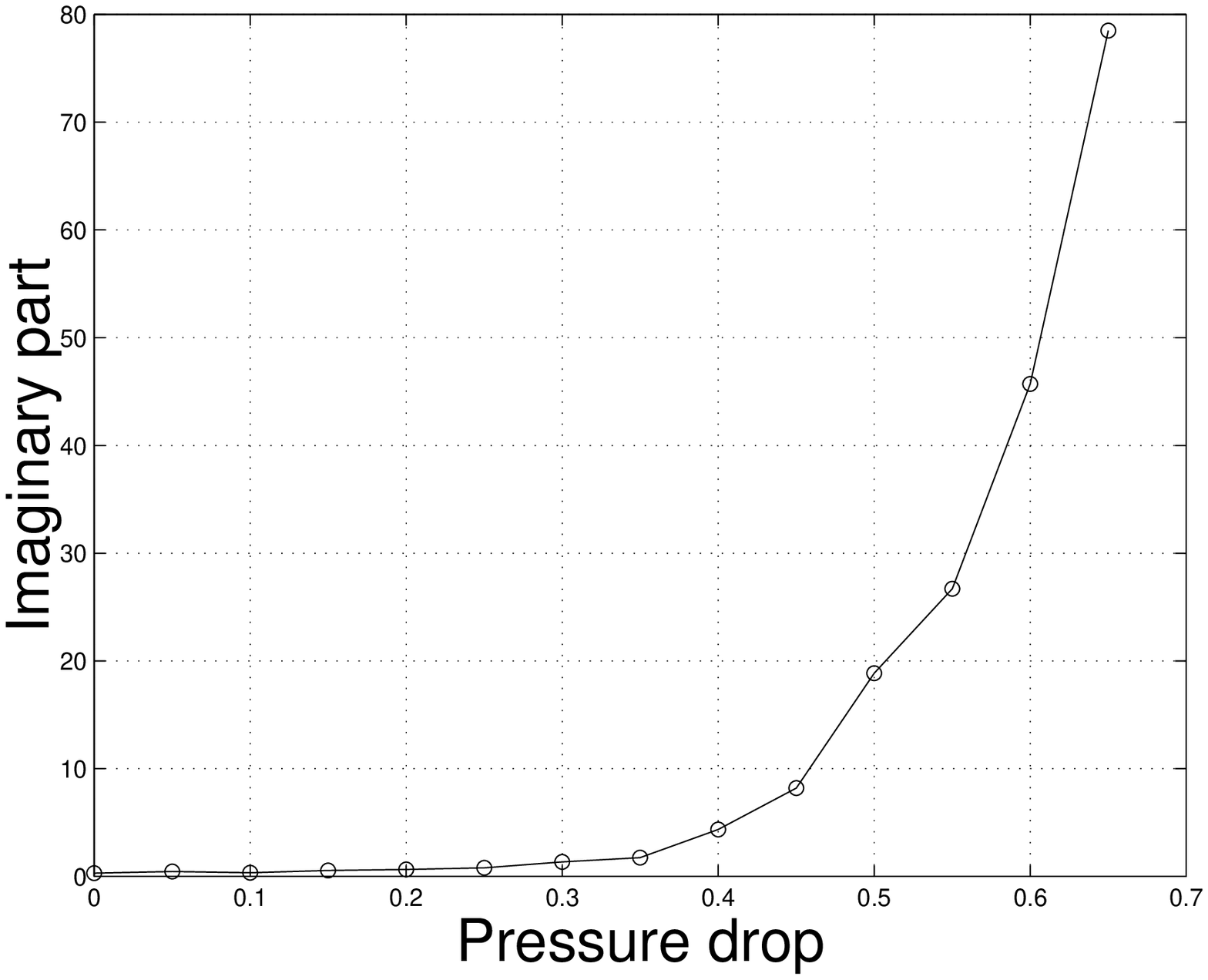}}

\begin{minipage}{10cm}
\caption{\it Real and imaginary parts of the critical eigenvalue as function of the control variable}
\end{minipage}
\end{center}
\end{figure}


\label{lastpage}


\begin{thebibliography}{aa}

\bibitem{AinsworthOden97}
    M.~Ainsworth and J.~T. Oden,
    \newblock A posteriori error estimation in finite element
    analysis,
    \newblock {\em Comput. Methods Appl. Mech. Engrg.}, 142:1--88, 1997.

\bibitem{BangerthRannacher02}
    W.~Bangerth and R.~Rannacher,
    \newblock Adaptive Methods for Differential Equations,
    \newblock Birkh\"auser, Basel, 2002, to appear.

\bibitem{Becker00}
    R.~Becker,
    \newblock An optimal-control approach to a posteriori error estimation
    for finite element discretizations of the Navier-Stokes
    equations,
    \newblock {\em East-West J. Numer. Math.}, 9:257--274, 2000.

\bibitem{Becker01}
    R.~Becker,
    \newblock Mesh adaptation for stationary flow control,
    \newblock {\em J. Math. Fluid Mech.}, 3:317--341, 2001.

\bibitem{BeckerRannacher95}
    R.~Becker and R.~Rannacher,
    \newblock Weighted a posteriori error control in {FE} methods,
    \newblock Lecture at ENUMATH-95, Paris, Sept. 18--22, 1995,
    Preprint 96-01, SFB 359, University of Heidelberg,
    \newblock Proc. {\em ENUMATH'97} (H.~G.~Bock {\it et al.}, eds),
    621--637, World Scientific, Singapore, 1998.

\bibitem{BeckerRannacher01}
    R.~Becker and R.~Rannacher,
    \newblock An optimal control approach to error estimation and mesh adaptation
    in finite element methods,
    \newblock {\em Acta Numerica 2000} (A. Iserles, ed.), 1--101, Cambridge
    University Press, 2001.

\bibitem{ErikssonEstepHansboJohnson95}
    K.~Eriksson, D.~Estep, P.~Hansbo, and C.~Johnson,
    \newblock Introduction to adaptive methods for differential equations.
    \newblock {\em Acta Numerica 1995} (A.~Iserles, ed.), 105--158,
    Cambridge University Press, 1995.

\bibitem{HeuvelineRannacher02}
    V.~Heuveline and R.~Rannacher,
    \newblock Adaptive finite element discretization of eigenvalue problems
    in hydrodynamic stability theory,
    \newblock Preprint, SFB 359, Universit\"at Heidelberg, March 2001.

\bibitem{Rannacher00} R.~Rannacher,
    {\em Finite element methods for the incompressible Navier-Stokes
    equations},
    Fundamental Directions in Mathematical Fluid Mechanics (G.~P. Galdi,
    et al., eds), 191--293, Birkh\"auser, Basel, 2000.

\bibitem{Verfurth96}
    R.~Verf\"urth,
    \newblock {\em A Review of A Posteriori Error Estimation and Adaptive
    Mesh-Refinement Techniques},
    \newblock Wiley/Teubner, New York Stuttgart, 1996.

\end{thebibliography}
\end{document}